\documentclass{article}
\usepackage{amssymb, amsmath, amsfonts}

\begin{document}

\begin{center}
{\Large Some remarks regarding finite bounded commutative BCK-algebras}

\begin{equation*}
\end{equation*}

\bigskip Cristina Flaut, \v{S}\'{a}rka Ho\v{s}kov\'{a}-Mayerov\'{a} and Radu
Vasile%
\begin{equation*}
\end{equation*}
\end{center}

\textbf{Abstract.} {\small \ In this chapter, starting from some results
obtained in the papers [FV; 19], [FHSV; 19], we provide some examples of
finite bounded commutative BCK-algebras, using the Wajsberg algebra
associated to a bounded commutative BCK-algebra. This method is an
alternative to the Iseki's construction, since by Iseki's extension some
properties of the obtained algebras are lost.}

\begin{equation*}
\end{equation*}%
\begin{equation*}
\end{equation*}

\textbf{1.} \textbf{Introduction}%
\begin{equation*}
\end{equation*}

BCK-algebras were first introduced in mathematics in 1966 by Y. Imai and K.
Iseki, through the paper [II; 66]. These algebras were presented as a
generalization \ of the concept of set-theoretic difference and
propositional calculi. The class of \ BCK-algebras is a proper subclass of
the class of BCI-algebras.\medskip

\textbf{Definition 1.1.} An algebra $(X,\ast ,\theta )$ of type $(2,0)$ is
called a \textit{BCI-algebra} if the following conditions are fulfilled:

$1)~((x\ast y)\ast (x\ast z))\ast (z\ast y)=\theta ,$ for all $x,y,z\in X;$

$2)~(x\ast (x\ast y))\ast y=\theta ,$ for all $x,y\in X;$

$3)~x\ast x=\theta ,$ for all $x\in X$;

$4)$ For all $x,y,z\in X$ such that $x\ast y=\theta ,y\ast x=\theta ,$ it
results $x=y$.

If a BCI-algebra $X$ satisfies the following identity:

$5)$ $\theta \ast x=\theta ,~$for all $x\in X,$ then $X$ is called a \textit{%
BCK-algebra}.\medskip

\textbf{Definition 1.2. }i) A BCK-algebra $(X,\ast ,\theta )$ is called 
\textit{commutative }if 
\begin{equation*}
x\ast (x\ast y)=y\ast (y\ast x),
\end{equation*}%
for all $x,y\in X$ and \textit{implicative }if 
\begin{equation*}
x\ast (y\ast x)=x,
\end{equation*}%
for all $x,y\in X.$

ii) ([Du; 99]) A BCK-algebra $(X,\ast ,\theta )$ is called \textit{positive
implicative} if and only if%
\begin{equation*}
\left( x\ast y\right) \ast z=\left( x\ast z\right) \ast (y\ast z),
\end{equation*}%
for all $x,y,z\in X.$

The \textit{partial order} relation on a BCK-algebra is defined such that $%
x\leq y$ if and only if $x\ast y=\theta .$

If in the BCK-algebra $(X,\ast ,\theta )$ there is an element $1$ such that $%
x\leq 1,$ for all $x\in X,$ then the algebra $X$ is called a \textit{bounded
BCK-algebra}. In a bounded BCK-algebra, we denote $1\ast x=\overline{x}.$

If in the bounded \ BCK-algebra $X$, an element $x\in X$ satisfies the
relation 
\begin{equation*}
\overline{\overline{x}}=x,
\end{equation*}%
then the element $x$ is called an \textit{involution}.

If $(X,\ast ,\theta )$ and $(Y,\circ ,\theta )$ are two BCK-algebras, a map $%
f:X\rightarrow Y$ with the property $f\left( x\ast y\right) =f\left(
x\right) \circ f\left( y\right) ,$ for all $x,y\in X,$ is called a \textit{%
BCK-algebras morphism}$.$ If $f$ is a bijective map, then $f$ is an \textit{%
isomorphism} of BCK-algebras.\medskip

\textbf{Definition 1.3.} 1) Let $(X,\ast ,\theta )$ be a BCK algebra and $Y$
be a nonempty subset of $X$. Therefore, $Y$ is called a \textit{subalgebra}
of the algebra $(X,\ast ,\theta )$ if and only if for each $x,y\in Y$, we
have $x\ast y\in Y$. This implies that $Y$ is closed to the binary
multiplication "$\ast $". \ 

It is well known that each BCK-algebra of degree $n+1$ contains a subalgebra
of degree $n.$

2) Let $(X,\ast ,\theta )$ be a BCK algebra and $I$ be a nonempty subset of $%
X$. Therefore, $I$ is called an \textit{ideal} of the algebra $X$ if and
only if for each $x,y\in X$ we have:

i) $\theta \in I;$

ii) $x\ast y\in I$ and $y\in I,$ then $x\in I$.\medskip\ 

\textbf{Proposition 1.4.} ([Me-Ju; 94]) \textit{Let} $(X,\ast ,\theta )$ 
\textit{be a BCK algebra and} $Y$ \textit{be a subalgebra of the algebra} $%
(X,\ast ,\theta )$. \textit{The following statements are true:}

\textit{i)} $\theta \in Y;$

\textit{ii)} $(Y,\ast ,\theta )$ \textit{is also a BCK-algebra.}$\Box
\medskip $

Let \thinspace $X$ be a BCK-algebra, such that $1\notin X$. On the set $%
Y=X\cup \{1\},$ we define the following multiplication $"\circ "$, as
follows:%
\begin{equation*}
x\circ y=\left\{ 
\begin{array}{c}
x\ast y,~if~x,y\in X; \\ 
\theta ,~if~x\in X,y=1; \\ 
1,~if~x=1~and~y\in X; \\ 
\theta ,~if~x=y=1.%
\end{array}%
\right.
\end{equation*}

The obtained algebra $\left( Y,\circ ,\theta \right) $ is a bounded
BCK-algebra and is obtained by the so-called Iseki's extension. The algebra $%
\left( Y,\circ ,\theta \right) $ is called t\textit{he algebra obtained from
algebra} $\left( X,\ast ,\theta \right) $ \textit{by Iseki's extension}
([Me-Ju; 94], Theorem 3.6).\medskip

\textbf{Remark 1.5.} ([Me-Ju; 94])

i) The Iseki's extension of a positive implicative BCK-algebra is still a
positive implicative BCK-algebra.

ii) The Iseki's \ extension of a commutative BCK-algebra, in general, is not
a commutative BCK-algebra.

iii) Let $X$ be a BCK-algebra and $Y$ its Iseki's extension. Therefore $X$
is an ideal in $Y$.

iv) If $I$ is an ideal of the BCK-algebra $X$, $x\in I$ and $y\leq x$, then $%
y\in I$.\medskip 

In the following, we will give some examples of finite bounded commutative
BCK-algebras. In the finite case, it is very useful to have many examples of
such algebras. But, such examples, in general, are not so easy to found. A
method for this purpose can be Iseki's extension. But, from the above, we
remark that the Iseki's extension can't be always used to obtain examples of
finite commutative bounded BCK-algebras with given initial properties, since
the commutativity, or other properties, can be lost. From this reason, we
use other technique to provide examples of such algebras. We use the
connections between finite commutative bounded BCK-algebras and Wajsberg
algebras and the algorithm and examples given in the papers [FHSV; 19] and
[FV; 19].\medskip 
\begin{equation*}
\end{equation*}

\textbf{2. Connections between finite bounded commutative BCK-algebras and
Wajsberg algebras}

\begin{equation*}
\end{equation*}

\textbf{Definition 2.1. \ }([CHA; 58]) An abelian monoid $\left( X,\theta
,\oplus \right) $ is called \textit{MV-algebra} if and only if we have an
unary operation $"^{\prime }"$ such that:

i) $(x^{\prime })^{\prime }=x;$

ii) $x\oplus \theta ^{\prime }=\theta ^{\prime };$

iii) $\left( x^{\prime }\oplus y\right) ^{\prime }\oplus y=$ $\left(
y^{\prime }\oplus x\right) ^{\prime }\oplus x$, for all $x,y\in X.$([Mu;
07]). We denote it by $\left( X,\oplus ,^{\prime },\theta \right) .\medskip $

We remark that in an MV-algebra the constant element $\theta ^{\prime }$ is
denoted with $1$. This is equivalent with

\begin{equation*}
1=\theta ^{\prime }.
\end{equation*}%
With the above definitions, the following multiplications are also defined:

\begin{equation*}
x\odot y=\left( x^{\prime }\oplus y^{\prime }\right) ^{\prime },
\end{equation*}%
\begin{equation*}
x\ominus y=x\odot y^{\prime }=\left( x^{\prime }\oplus y\right) ^{\prime }.
\end{equation*}%
(see ([Mu; 07]))

From [COM; 00], Theorem 1.7.1, for a bounded commutative BCK-algebra $\left(
X,\ast ,\theta ,1\right) $, if we define 
\begin{equation*}
x^{\prime }=1\ast x,
\end{equation*}%
\begin{equation*}
x\oplus y=1\ast \left( \left( 1\ast x\right) \ast y\right) =\left( x^{\prime
}\ast y\right) ^{\prime },x,y\in X,
\end{equation*}%
we obtain that the algebra $\left( X,\oplus ,^{\prime },\theta \right) $ is
an \textit{MV}-algebra, with 
\begin{equation*}
x\ominus y=x\ast y.
\end{equation*}

The converse is also true, that means if $\left( X,\oplus ,\theta ,^{\prime
}\right) $ is an MV-algebra, then $\left( X,\ominus ,\theta ,1\right) $ is a
bounded commutative BCK-algebra.\medskip

\textbf{Definition 2.2. }([COM; 00], Definition 4.2.1) An algebra $\left(
W,\circ ,\overline{\phantom{x}}~,1\right) $ of type $\left( 2,1,0\right) ~$%
is called a \textit{Wajsberg algebra (}or\textit{\ W-algebra)} if and only
if the following conditions are fulfilled:

i) $1\circ x=x;$

ii) $\left( x\circ y\right) \circ \left[ \left( y\circ z\right) \circ \left(
x\circ z\right) \right] =1;$

iii) $\left( x\circ y\right) \circ y=\left( y\circ x\right) \circ x;$

iv) $\left( \overline{x}\circ \overline{y}\right) \circ \left( y\circ
x\right) =1$, for every $x,y,z\in W$.\medskip

\textbf{Remark 2.3. }([COM; 00], Lemma 4.2.2 and Theorem 4.2.5)

i) For the Wajsberg algebra $\left( W,\circ ,\overline{\phantom{x}},1\right) 
$, if we define the following multiplications 
\begin{equation*}
x\odot y=\overline{\left( x\circ \overline{y}\right) }
\end{equation*}%
and 
\begin{equation*}
x\oplus y=\overline{x}\circ y,
\end{equation*}%
for all $x,y\in W$, we obtain that $\left( W,\oplus ,\odot ,\overline{%
\phantom{x}},\theta ,1\right) $ is an MV-algebra.

ii) Conversely, if $\left( X,\oplus ,\odot ,^{\prime },\theta ,1\right) $ is
an MV-algebra, defining on $X$ the operation%
\begin{equation*}
x\circ y=x^{\prime }\oplus y,
\end{equation*}%
we obtain that $\left( X,\circ ,^{\prime },1\right) $ is a Wajsberg
algebra.\medskip

\textbf{Remark 2.4}. From the above, if\ $\left( W,\circ ,\overline{%
\phantom{x}},1\right) $ is a Wajsberg algebra, then $\left( W,\oplus ,\odot ,%
\overline{\phantom{x}},0,1\right) $ is an MV-algebra, with 
\begin{equation}
x\ominus y=\overline{\left( \overline{x}\oplus y\right) }=\overline{\left(
x\circ y\right) }.  \tag{2.1.}
\end{equation}%
Defining 
\begin{equation}
x\ast y=\overline{\left( x\circ y\right) },  \tag{2.2.}
\end{equation}%
we have that $(W,\ast ,\theta ,1)$ is a bounded commutative
BCK-algebra.\medskip 

Using the above remark, starting from some known finite examples of Wajsberg
algebras given in the papers [FHSV; 19] and [FV; 19], we can obtain examples
of finite commutative bounded BCK-algebras, using the following
algorithm.\medskip

\textbf{The Algorithm\medskip }

1) Let $n$ be a natural number, $n\neq 0$ and 
\begin{equation*}
n=r_{1}r_{2}...r_{t},r_{i}\in \mathbb{N},1<r_{i}<n,i\in \{1,2,...,t\},
\end{equation*}%
be the decomposition of the number $n$ in factors. The decompositions with
the same terms, but with other order of them in the product, will be counted
one time. The number of all such decompositions will be denoted with $\pi
_{n}$.

2) There are only $\pi _{n}$ nonismorphic, as ordered sets, Wajsberg
algebras with $n$ elements. We obtain these algebras as a finite product of
totally ordered Wajsberg algebras (see [BV; 10] and \textbf{[}FV; 19\textbf{%
], }Theorem 4.8 \textbf{)}.

3) Using Remark 2.4 from above, a commutative bounded BCK-algebra can be
associated to each Wajsberg algebra.

\begin{equation*}
\end{equation*}

\textbf{3. Examples of finite commutative bounded BCK-algebras}

\begin{equation*}
\end{equation*}

In the following, we use some examples of Wajsberg algebras given in the
paper [FHSV; 19]. To these algebras, we will associate the corresponding
commutative bounded BCK-algebras and we give their subalgebras and
ideals.\medskip

\textbf{Example 3.1.} Let $W=\{O\leq A\leq B\leq E\}$ be a totally ordered
set. On $W$ we define the multiplication $\circ _{1}~$as in the below table,
such that $\left( W,\circ _{1},E\right) $ is a Wajsberg algebra. We have $%
\overline{A}=B$ and $\overline{B}=A$.

\begin{equation*}
\begin{tabular}{l|llll}
$\circ _{1}$ & $O$ & $A$ & $B$ & $E$ \\ \hline
$O$ & $E$ & $E$ & $E$ & $E$ \\ 
$A$ & $B$ & $E$ & $E$ & $E$ \\ 
$B$ & $A$ & $B$ & $E$ & $E$ \\ 
$E$ & $O$ & $A$ & $B$ & $E$%
\end{tabular}%
\text{.}
\end{equation*}%
(see [FHSV; 19], Example 4.1.1)

Therefore, the associated commutative bounded BCK-algebras $\left( W,\ast
_{1},O\right) ~$has multiplication given in the below table:%
\begin{equation*}
\begin{tabular}{l|llll}
$\ast _{1}$ & $O$ & $A$ & $B$ & $E$ \\ \hline
$O$ & $O$ & $O$ & $O$ & $O$ \\ 
$A$ & $A$ & $O$ & $O$ & $O$ \\ 
$B$ & $B$ & $A$ & $O$ & $O$ \\ 
$E$ & $E$ & $B$ & $A$ & $O$%
\end{tabular}%
\text{.}
\end{equation*}

The proper subalgebras of this algebra are: $\{O,A\},\{O,B\},\{O,E\},\{O,A,B%
\}$. There are no proper ideals in the algebra $\left( W,\ast _{1},O\right) $%
. \medskip

\textbf{Example 3.2.} Let $W=\{O\leq A\leq B\leq E\}$ be a totally ordered
set. On $W$ we define the multiplication $\circ _{2}$ as in the below table,
such that $\left( W,\circ _{2},E\right) $ is a Wajsberg algebra. 
\begin{equation*}
\begin{tabular}{l|llll}
$\circ _{2}$ & $O$ & $A$ & $B$ & $E$ \\ \hline
$O$ & $E$ & $E$ & $E$ & $E$ \\ 
$A$ & $B$ & $E$ & $B$ & $E$ \\ 
$B$ & $A$ & $A$ & $E$ & $E$ \\ 
$E$ & $O$ & $A$ & $B$ & $E$%
\end{tabular}%
\text{.}
\end{equation*}

Therefore, the associated commutative bounded BCK-algebras $\left( W,\ast
_{2},O\right) $ has multiplication given below: 
\begin{equation*}
\begin{tabular}{l|llll}
$\ast _{2}$ & $O$ & $A$ & $B$ & $E$ \\ \hline
$O$ & $O$ & $O$ & $O$ & $O$ \\ 
$A$ & $A$ & $O$ & $A$ & $O$ \\ 
$B$ & $B$ & $B$ & $O$ & $O$ \\ 
$E$ & $E$ & $B$ & $A$ & $O$%
\end{tabular}%
\text{.}
\end{equation*}%
The proper subalgebras of this algebra are: $\{O,A\},\{O,B\},\{O,E\},\{O,A,B%
\}$. The proper ideals are: $\{O,A\},\{O,B\}$.\medskip

\textbf{Example 3.3.} Let $W=\{O\leq A\leq B\leq C\leq D\leq E\}$ be a
totally ordered set. On $W$ we define a multiplication $\circ _{3}~$given in
the below table, such that $\left( W,\circ _{3},E\right) $ is a Wajsberg
algebra. We have $\overline{A}=D$, $\overline{B}=C$, $\overline{C}=B$, $%
\overline{D}=A$. 
\begin{equation*}
\begin{tabular}{l|llllll}
$\circ _{3}$ & $O$ & $A$ & $B$ & $C$ & $D$ & $E$ \\ \hline
$O$ & $E$ & $E$ & $E$ & $E$ & $E$ & $E$ \\ 
$A$ & $D$ & $E$ & $E$ & $E$ & $E$ & $E$ \\ 
$B$ & $C$ & $D$ & $E$ & $E$ & $E$ & $E$ \\ 
$C$ & $B$ & $C$ & $D$ & $E$ & $E$ & $E$ \\ 
$D$ & $A$ & $B$ & $C$ & $D$ & $E$ & $E$ \\ 
$E$ & $O$ & $A$ & $B$ & $C$ & $D$ & $E$%
\end{tabular}%
.
\end{equation*}%
Therefore, the associated commutative bounded BCK-algebras $\left( W,\ast
_{3},O\right) $ has multiplication defined below:

\begin{equation*}
\begin{tabular}{l|llllll}
$\ast _{3}$ & $O$ & $A$ & $B$ & $C$ & $D$ & $E$ \\ \hline
$O$ & $O$ & $O$ & $O$ & $O$ & $O$ & $O$ \\ 
$A$ & $A$ & $O$ & $O$ & $O$ & $O$ & $O$ \\ 
$B$ & $B$ & $A$ & $O$ & $O$ & $O$ & $O$ \\ 
$C$ & $C$ & $B$ & $A$ & $O$ & $O$ & $O$ \\ 
$D$ & $D$ & $C$ & $B$ & $A$ & $O$ & $O$ \\ 
$E$ & $E$ & $D$ & $C$ & $B$ & $A$ & $O$%
\end{tabular}%
.
\end{equation*}%
The proper subalgebras of this algebra are: $\{O,A\},\{O,B\},\{O,C\},\{O,D%
\}, $\newline
$\{O,E\},\{O,A,B\},\{O,B,D\},\{O,A,B,C\},\{O,A,B,C,D\}$ . \newline
This algebra has no proper ideals.\medskip

\textbf{Example 3.4. }Let $W=\{O\leq A\leq B\leq C\leq D\leq E\}$ be a
totally ordered set. On $W$ we define a multiplication $\circ _{4}~$given in
the below table, such that $\left( W,\circ _{4},E\right) $ is a Wajsberg
algebra. We have 
\begin{equation*}
\begin{tabular}{l|llllll}
$\circ _{4}$ & $O$ & $A$ & $B$ & $C$ & $D$ & $E$ \\ \hline
$O$ & $E$ & $E$ & $E$ & $E$ & $E$ & $E$ \\ 
$A$ & $D$ & $E$ & $E$ & $D$ & $E$ & $E$ \\ 
$B$ & $C$ & $D$ & $E$ & $C$ & $D$ & $E$ \\ 
$C$ & $B$ & $B$ & $B$ & $E$ & $E$ & $E$ \\ 
$D$ & $A$ & $B$ & $B$ & $D$ & $E$ & $E$ \\ 
$E$ & $O$ & $A$ & $B$ & $C$ & $D$ & $E$%
\end{tabular}%
\ .
\end{equation*}%
Therefore, the associated commutative bounded BCK-algebras $\left( W,\ast
_{4},O\right) $ has multiplication given in the following table: 
\begin{equation*}
\begin{tabular}{l|llllll}
$\ast _{4}$ & $O$ & $A$ & $B$ & $C$ & $D$ & $E$ \\ \hline
$O$ & $O$ & $O$ & $O$ & $O$ & $O$ & $O$ \\ 
$A$ & $A$ & $O$ & $O$ & $A$ & $O$ & $O$ \\ 
$B$ & $B$ & $A$ & $O$ & $B$ & $A$ & $O$ \\ 
$C$ & $C$ & $C$ & $C$ & $O$ & $O$ & $O$ \\ 
$D$ & $D$ & $C$ & $C$ & $A$ & $O$ & $O$ \\ 
$E$ & $E$ & $D$ & $C$ & $B$ & $A$ & $O$%
\end{tabular}%
\ \text{.}
\end{equation*}%
The proper subalgebras of this algebra are: $\{O,A\},\{O,B\},\{O,C\},\{O,D%
\}, $\newline
$\{O,E\},\{O,A,B\},\{O,A,B,C\},\{O,A,B,C,D\},\{O,A,C\},\{O,A,C,D\}$. \newline
The proper ideals of this algebra are:\medskip\ $\{O,A,B\},\{O,C\}$.\medskip

\textbf{Example 3.5. }Let $W=\{O\leq A\leq B\leq C\leq D\leq E\}$ be a
totally ordered set. On $W$ we define a multiplication $\circ _{5}~$given in
the below table, such that $\left( W,\circ _{5},E\right) $ is a Wajsberg
algebra. We have 
\begin{equation*}
\begin{tabular}{l|llllll}
$\circ _{5}$ & $O$ & $A$ & $B$ & $C$ & $D$ & $E$ \\ \hline
$O$ & $E$ & $E$ & $E$ & $E$ & $E$ & $E$ \\ 
$A$ & $C$ & $E$ & $A$ & $D$ & $D$ & $E$ \\ 
$B$ & $D$ & $E$ & $E$ & $D$ & $D$ & $E$ \\ 
$C$ & $A$ & $E$ & $A$ & $E$ & $E$ & $E$ \\ 
$D$ & $B$ & $A$ & $B$ & $A$ & $E$ & $E$ \\ 
$E$ & $O$ & $A$ & $B$ & $C$ & $D$ & $E$%
\end{tabular}%
.
\end{equation*}%
Therefore, the associated commutative bounded BCK-algebras $\left( W,\ast
_{5},O\right) $ has multiplication defined in the following table:%
\begin{equation*}
\begin{tabular}{l|llllll}
$\ast _{5}$ & $O$ & $A$ & $B$ & $C$ & $D$ & $E$ \\ \hline
$O$ & $O$ & $O$ & $O$ & $O$ & $O$ & $O$ \\ 
$A$ & $A$ & $O$ & $C$ & $B$ & $B$ & $O$ \\ 
$B$ & $B$ & $O$ & $O$ & $B$ & $B$ & $O$ \\ 
$C$ & $C$ & $O$ & $C$ & $O$ & $O$ & $O$ \\ 
$D$ & $D$ & $C$ & $D$ & $C$ & $O$ & $O$ \\ 
$E$ & $E$ & $C$ & $D$ & $A$ & $B$ & $O$%
\end{tabular}%
.
\end{equation*}

The proper subalgebras of this algebra are: $\{O,A\},\{O,B\},\{O,C\},\{O,D%
\}, $\newline
$\{O,E\},\{O,B,C\},\{O,C,D\},\{O,A,B,C\},\{O,A,B,C,D\}.$

All proper ideals are: $\{O,C,D\}$, $\{O,B\}$.\medskip

\textbf{Example 3.6. }Let $W=\{O\leq X\leq Y\leq Z\leq T\leq U\leq V\leq E\}$
be a totally ordered set. On $W$ we define a multiplication $\circ _{6}~$%
which can be found in the below table. We obtain that $\left( W,\circ
_{6},E\right) $ is a Wajsberg algebra.We have $\overline{X}=V$, $\overline{Y}%
=U$, $\overline{Z}=T$. Therefore the algebra $W$ has the following
multiplication table:%
\begin{equation*}
\begin{tabular}{l|llllllll}
$\circ _{6}$ & $O$ & $X$ & $Y$ & $Z$ & $T$ & $U$ & $V$ & $E$ \\ \hline
$O$ & $E$ & $E$ & $E$ & $E$ & $E$ & $E$ & $E$ & $E$ \\ 
$X$ & $V$ & $E$ & $E$ & $E$ & $E$ & $E$ & $E$ & $E$ \\ 
$Y$ & $U$ & $V$ & $E$ & $E$ & $E$ & $E$ & $E$ & $E$ \\ 
$Z$ & $T$ & $U$ & $V$ & $E$ & $E$ & $E$ & $E$ & $E$ \\ 
$T$ & $Z$ & $T$ & $U$ & $V$ & $E$ & $E$ & $E$ & $E$ \\ 
$U$ & $Y$ & $Z$ & $T$ & $U$ & $V$ & $E$ & $E$ & $E$ \\ 
$V$ & $X$ & $Y$ & $Z$ & $T$ & $U$ & $V$ & $E$ & $E$ \\ 
$E$ & $O$ & $X$ & $Y$ & $Z$ & $T$ & $U$ & $V$ & $E$%
\end{tabular}%
.
\end{equation*}%
From here, we get thet the associated commutative bounded BCK-algebras $%
\left( W,\ast _{6},O\right) $ has multiplication given in the below table:%
\begin{equation*}
\begin{tabular}{l|llllllll}
$\ast _{6}$ & $O$ & $X$ & $Y$ & $Z$ & $T$ & $U$ & $V$ & $E$ \\ \hline
$O$ & $O$ & $O$ & $O$ & $O$ & $O$ & $O$ & $O$ & $O$ \\ 
$X$ & $X$ & $O$ & $O$ & $O$ & $O$ & $O$ & $O$ & $O$ \\ 
$Y$ & $Y$ & $X$ & $O$ & $O$ & $O$ & $O$ & $O$ & $O$ \\ 
$Z$ & $Z$ & $Y$ & $X$ & $O$ & $O$ & $O$ & $O$ & $O$ \\ 
$T$ & $T$ & $Z$ & $Y$ & $X$ & $O$ & $O$ & $O$ & $O$ \\ 
$U$ & $U$ & $T$ & $Z$ & $Y$ & $X$ & $O$ & $O$ & $O$ \\ 
$V$ & $V$ & $U$ & $T$ & $Z$ & $Y$ & $X$ & $O$ & $O$ \\ 
$E$ & $E$ & $V$ & $U$ & $T$ & $Z$ & $T$ & $X$ & $O$%
\end{tabular}%
\text{.}
\end{equation*}%
The proper subalgebras of this algebra are: $\{O,J\},J\in \{X,Y,Z,T,U,V,E\},$%
\newline
$\{O,X,Y\},\{O,X,Y,Z\},\{O,X,Y,Z,T\},\{O,X,Y,Z,T,U\},\{O,X,Y,Z,T,U,V\}$.
There are no proper ideals.\medskip

\textbf{Example 3.7. }Let $W=\{O\leq X\leq Y\leq Z\leq T\leq U\leq V\leq E\}$
be a totally ordered set. On $W$ we define a multiplication $\circ _{7}~$%
given in the below table, such that $\left( W,\circ _{7},E\right) $ is a
Wajsberg algebra. \textbf{\ } 
\begin{equation*}
\begin{tabular}{l|llllllll}
$\circ _{7}$ & $O$ & $X$ & $Y$ & $Z$ & $T$ & $U$ & $V$ & $E$ \\ \hline
$O$ & $E$ & $E$ & $E$ & $E$ & $E$ & $E$ & $E$ & $E$ \\ 
$X$ & $V$ & $E$ & $V$ & $E$ & $V$ & $E$ & $V$ & $E$ \\ 
$Y$ & $U$ & $U$ & $E$ & $E$ & $E$ & $E$ & $E$ & $E$ \\ 
$Z$ & $T$ & $U$ & $V$ & $E$ & $V$ & $E$ & $V$ & $E$ \\ 
$T$ & $Z$ & $Z$ & $U$ & $U$ & $E$ & $E$ & $E$ & $E$ \\ 
$U$ & $Y$ & $Z$ & $T$ & $U$ & $T$ & $E$ & $V$ & $E$ \\ 
$V$ & $X$ & $X$ & $Z$ & $Z$ & $U$ & $U$ & $E$ & $E$ \\ 
$E$ & $O$ & $X$ & $Y$ & $Z$ & $T$ & $U$ & $V$ & $E$%
\end{tabular}%
.
\end{equation*}%
Therefore, the associated commutative bounded BCK-algebras $\left( W,\ast
_{7},O\right) $ has multiplication defined in the below table:%
\begin{equation*}
\begin{tabular}{l|llllllll}
$\ast _{7}$ & $O$ & $X$ & $Y$ & $Z$ & $T$ & $U$ & $V$ & $E$ \\ \hline
$O$ & $O$ & $O$ & $O$ & $O$ & $O$ & $O$ & $O$ & $O$ \\ 
$X$ & $X$ & $O$ & $X$ & $O$ & $X$ & $O$ & $X$ & $O$ \\ 
$Y$ & $Y$ & $Y$ & $O$ & $O$ & $O$ & $O$ & $O$ & $O$ \\ 
$Z$ & $Z$ & $Y$ & $X$ & $O$ & $X$ & $O$ & $X$ & $O$ \\ 
$T$ & $T$ & $T$ & $Y$ & $Y$ & $O$ & $O$ & $O$ & $O$ \\ 
$U$ & $U$ & $T$ & $Z$ & $Y$ & $Z$ & $O$ & $X$ & $O$ \\ 
$V$ & $V$ & $V$ & $T$ & $T$ & $Y$ & $Y$ & $O$ & $O$ \\ 
$E$ & $E$ & $V$ & $U$ & $T$ & $Z$ & $Y$ & $X$ & $O$%
\end{tabular}%
.
\end{equation*}%
The proper subalgebras of this algebra are: $\{O,J\},J\in \{X,Y,Z,T,U,V,E\},$%
\newline
$\{O,X,Y\},\{O,X,Y,Z\},\{O,T,Y\},\{O,T,Y,V\},\{O,X,Y,Z,T\},\{O,X,Y,Z,T,U\},$%
\newline
$\{O,X,Y,Z,T,U,V\}.$

All proper ideals are:$~\{O,Y,T,V\}$, $\{O,X\}$.\medskip

\textbf{Conclusions.} In this chapter, we provided an algorithm for finding
examples of \ finite commutative bounded BCK-algebras, using their
connections with Wajsberg algebras. This algorithm allows us to find such
examples no matter the order of the algebra. This thing is very useful,
since examples of such algebras are very rarely encountered in the specialty
books.

\begin{equation*}
\end{equation*}

\bigskip \textbf{References}%
\begin{equation*}
\end{equation*}

[AAT; 96] Abujabal, H.A.S., Aslam, M., Thaheem, A.B., \textit{A
representation of bounded commutative BCK-algebras}, Internat. J. Math. \&
Math. Sci., 19(4)(1996), 733-736.

[BV; 10] Belohlavek, R., Vilem Vychodil, V., \textit{Residuated Lattices of
Size} $\leq 12$, Order, 27(2010), 147-161.

[CHA; 58] Chang, C.C.,\textit{\ Algebraic analysis of many-valued logic},
Trans. Amer. Math. Soc. 88(1958), 467-490.

[COM; 00] Cignoli, R. L. O, Ottaviano, I. M. L. D, Mundici, D., \textit{%
Algebraic foundations of many-valued reasoning}, Trends in Logic, Studia
Logica Library, Dordrecht, Kluwer Academic Publishers, 7(2000).

[Du; 99] W.A. Dudek, \ \textit{On embedding Hilbert algebras in BCK-algebras}%
, Mathematica Moravica, \textbf{3(1999)}, 25-28.

[FHSV; 19] Flaut, C., Ho\v{s}kov\'{a}-Mayerov\'{a}, \v{S}., Saeid, A., B.,
Vasile, R., \textit{Wajsberg algebras of order} $n,n\leq 9$,
https://arxiv.org/pdf/1905.05755.pdf

[FV; 19] Flaut, C., Vasile, R., \textit{Wajsberg algebras arising from
binary block codes}, https://arxiv.org/pdf/1904.07169.pdf

[II; 66] Imai, Y., Iseki, K., \textit{On axiom systems of propositional
calculi}, Proc. Japan Academic, \textbf{42(1966)}, 19-22.

[JS; 11] Jun, Y. B., Song, S. Z., \textit{Codes based on BCK-algebras},
Inform. Sciences., \textbf{181(2011)}, 5102-5109.

[Me-Ju; 94] Meng, J., Jun, Y. B., \textit{BCK-algebras}, Kyung Moon Sa Co.
Seoul, Korea, 1994.

[Mu; 07] Mundici, D., \textit{MV-algebras-a short tutorial}, Department of
Mathematics \textit{Ulisse Dini}, University of Florence, 2007.%
\begin{equation*}
\end{equation*}

Cristina FLAUT

{\small Faculty of Mathematics and Computer Science, }\newline

{\small Ovidius University of Constan\c{t}a, Rom\^{a}nia,}

{\small Bd. Mamaia 124, 900527,}

{\small http://www.univ-ovidius.ro/math/}

{\small e-mail: cflaut@univ-ovidius.ro; cristina\_flaut@yahoo.com}%
\begin{equation*}
\end{equation*}

\v{S}\'{a}rka Ho\v{s}kov\'{a}-Mayerov\'{a}

{\small Department of Mathematics and Physics,}

{\small University of Defence, Brno, Czech Republic}

{\small e-mail: sarka.mayerova@unob.cz}%
\begin{equation*}
\end{equation*}

\qquad \qquad

Radu Vasile,

{\small PhD student at Doctoral School of Mathematics,}

{\small Ovidius University of Constan\c{t}a, Rom\^{a}nia}

{\small rvasile@gmail.com}

\end{document}